\documentclass[12pt,oneside,reqno]{amsart}

\usepackage[utf8]{inputenc}
\usepackage[a4paper,nomarginpar]{geometry}

\usepackage[english]{babel}
\usepackage{enumitem}
\usepackage{amsmath,amsfonts,amssymb,amsthm}
\usepackage{amscd}
\usepackage{tikz}

\allowdisplaybreaks

\makeatletter
\def\subsection{\@startsection{subsection}{2}%
	\z@{.5\linespacing\@plus.7\linespacing}{.3\linespacing}%
	{\normalfont\bfseries}}
\makeatother

\newtheorem{theorem}{Theorem}%[section]
\newtheorem{lemma}[theorem]{Lemma}
\newtheorem{proposition}[theorem]{Proposition}
\newtheorem{corollary}[theorem]{Corollary}

\theoremstyle{definition}

\begin{document}
	
	% \title[short text for running head]{full title}
	
	\title[On the entropy of processes generated by quasifactors]
	{On the entropy of processes generated by quasifactors}
	
	% Only \author and \address are required; other information is optional.
	% Remove any unused author tags.
	
	% author one information
	% \author[short version for running head]{name for top of paper}

    \author{R\^omulo M. Vermersch}
	\address{Departamento de Matem\'atica, Centro de Ci\^encias F\'isicas e
		Matem\'aticas, Universidade Federal de Santa Catarina,
		Florian\'opolis, SC, 88040-900, Brazil.}
	\curraddr{}
	\email{romulo.vermersch@ufsc.br}
	\thanks{}

	% author two information
	
%	\author{}
%	\address{}
%	\curraddr{}
%	\email{}
%	\thanks{}

	% \subjclass is required.
	
	
	\date{}
	
	\dedicatory{}
	\begin{abstract}
		
		Let $(X,\mathcal{X},\mu,T)$ be a measurable dynamical system, where $T$ is a homeomorphism acting on the zero-dimensional compact metric space $X$, $\mathcal{X}$ is the Borel $\sigma$-algebra on $X$ and $\mu$ is an ergodic $T$-invariant Borel probability measure on $X$. In this work we show that, if $h_{\mu}(T)>0$, then $h_{\widetilde{\mu}}(\widetilde{T})>0$ for every ergodic quasifactor $\widetilde{\mu}$ of $\mu$ having full-support. This result can be seen, in the zero-dimensional setting, as a counterpart of the Glasner-Weiss result about quasifactors of zero-entropy systems.
	\end{abstract}
    
	\maketitle

	%%%%%%%%%%%%%%%%%%%%%%%%%%%%%%%%%%%%%%%%%%%%%%%%%%%%%%%%%%%%%%%%%%%%%%%%%%%%%

	\section{Introduction}\label{Intro}

    Entropy is a central notion in ergodic theory, providing a fundamental measure of the unpredictability and complexity of a dynamical system. Actually, since Kolmogorov's foundational work \cite{Kolmogorov}, entropy has become a major theme within important directions of research such as: isomorphism theory \cite{Ornstein75}, Lyapunov exponents \cite{Katok80,Pesin77,Young}, volume growth rates \cite{Yomdin}, uniformly \cite{Anosov67,Smale} and non-uniformly \cite{Pesin77} hyperbolic dynamical systems. A detailed account of the deep connections between entropy and these topics can be found in the expository paper by Katok \cite{Katok07}, where many more references can also be found. 
    
    \smallskip
    
    By a {\em measurable dynamical system} (MDS) we mean a quadriple $\mathfrak{X}=(X,\mathcal{X},\mu,T)$, where $X$ is a compact metric space, $\mathcal{X}$ is the Borel $\sigma$-algebra on $X$, $\mu$ is a Borel probability measure on $\mathcal{X}$ and $T:X\to X$ is a homeomorphism that preserves $\mu$.
    
    By a {\em topological dynamical system} (TDS) we mean a pair $(X,T)$
    consisting of a compact metric space $X$ and a homeomorphism $T : X \to X$. 
    Such a TDS induces, in a natural way, the TDS $(\mathcal{M}(X),\widetilde{T}).$
    Here, $\mathcal{M}(X)$ denotes the space of all Borel probability
    measures on $X$ endowed with the {\em Prokhorov metric}
    $$
    d_P(\mu,\nu):= \inf\{\delta > 0 : \mu(A) \leqslant \nu(A^\delta) + \delta
    \text{ for all } A \in \mathcal{X}\},
    $$
    and
    $\widetilde{T} : \mathcal{M}(X) \to \mathcal{M}(X)$ is the homeomorphism given by
    $$
    (\widetilde{T}(\mu))(A):= \mu(T^{-1}(A)) \ \ \ (\mu \in \mathcal{M}(X), A \in \mathcal{X}).
    $$
    
    It is well known that  $\mathcal{M}(X)$ is a compact metric space and that $d_P(\mu,\nu)$ induces the so-called
    {\em weak*-topology} on $\mathcal{M}(X)$, that is, the topology whose basic open neighborhoods
    of $\mu \in \mathcal{M}(X)$ are the sets of the form
    $$
    \mathbb{V}(\mu;f_1,\ldots,f_k;\varepsilon):= \Big\{\nu \in \mathcal{M}(X) :
    \Big| \int_X f_i \,d\nu - \int_X f_i \,d\mu \Big| < \varepsilon
    \text{ for } i = 1,\ldots,k\Big\},
    $$
    where $k \geqslant 1$, $f_1,\ldots,f_k : X \to \mathbb{R}$ are continuous functions and
    $\varepsilon > 0$.
    
    We refer the reader to the books
    \cite{PBil99,RDud02,AKec95} for a study of the space $\mathcal{M}(X)$.
    
    \smallskip

    The research on the connections between the dynamics of the TDS $(X,T)$ and
    the dynamics of the induced TDS $(\mathcal{M}(X),\widetilde{T})$ was initiated by Bauer and Sigmund \cite{WBauKSig75},
    and was later developed by several authors; see \cite{NBerUDarRVer22,NBerRVer16,Burguet_Shi,EGlaBWei95,EGlaBWei03,LiOprWu17,LiYanYe15,QiaoZhou17,ShiZhang,KSig78}, for instance.
    The TDS $(\mathcal{M}(X),\widetilde{T})$ serves as an abstract model for systems in statistical mechanics, where the dynamics can be described in deterministic terms (in the sense that the time-evolution of the system is given by some physical law), but the states of the system are probability distributions on the phase space. 
    In \cite{EGla83} Glasner introduced the notion of a {\em quasifactor} of a MDS as an ergodic-theoretic analogue of an induced TDS. Let us see its definition:
    
    A {\it quasifactor} of $\mathfrak{X}=(X,\mathcal{X},\mu,T)$ is a MDS $\widetilde{\mathfrak{X}}=(\mathcal{M}(X),\widetilde{\mathcal{X}},\widetilde{\mu},\widetilde{T})$ such that  $\widetilde{\mu}$ satisfies the so-called {\it barycenter equation}:
    \begin{equation}\label{barycenter} 
    	\mu=\displaystyle\int_{\mathcal{M}(X)}\theta d\widetilde{\mu}(\theta)
    \end{equation}
    Here, $\widetilde{\mathcal{X}}$ denotes the Borel $\sigma$-algebra on $\mathcal{M}(X)$.
    Equivalently, we say that {\em $\mu$ is the barycenter of $\widetilde{\mu}$}.
    
    The barycenter equation means that, by choosing any compact topology on $X$ compatible with its Borel structure one has
    \begin{equation*}\displaystyle\int_{X}f(x)d \mu(x)=\int_{\mathcal{M}(X)}\int_{X}f(x)d \theta(x)d \widetilde{\mu}(\theta)
    \end{equation*}
    for all $f:X\rightarrow\mathbb{R}$ continuous function.
    
    Glasner also showed that this definition is independent of the choice of the compact topology compatible with the Borel structure (\cite{EGla83}). 
    For convenience, sometimes we will say that $\widetilde{\mu}$ is a quasifactor of $\mu$ and we shall denote by $Q(\mu)$ the set of all quasifactors of $\mu$ and by $Q_{e}(\mu)$ the set of all ergodic quasifactors of $\mu$.
    
    Moreover we remark that, for each fixed $A\in\mathcal{X}$ the map $\Psi_{A}:\nu\in\mathcal{M}(X)\mapsto \nu(A)\in [0,1]$ is Borel and $\displaystyle \mu(A)=\int_{\mathcal{M}(X)}\nu(A)d \widetilde{\mu}(\nu)$; for a proof of this well-known fact see Lemma 4.1 from \cite{HanfengLi_KairanLiu}.
    
    \smallskip
    
    In this work we are concerned with the relationship between the entropy of the MDS $\mathfrak{X}=(X,\mathcal{X},\mu,T)$ and of the MDS $\widetilde{\mathfrak{X}}=(\mathcal{M}(X),\widetilde{\mathcal{X}},\widetilde{\mu},\widetilde{T})$, where $\mu$ is ergodic and $\widetilde{\mu}\in Q_{e}(\mu)$ is an ergodic quasifactor of $\mu$.

    The research on the relationship between the entropy of a MDS and of a quasifactor of it can be traced back to a deep result due to Glasner and Weiss \cite{EGlaBWei95} which asserts that
    if $\mathfrak{X}=(X,\mathcal{X},\mu,T)$ has zero entropy, then so does $\widetilde{\mathfrak{X}}=(\mathcal{M}(X),\widetilde{\mathcal{X}},\widetilde{\mu},\widetilde{T})$ for {\em every} $\widetilde{\mu}\in Q(\mu)$. By the variational principle it implies that, if $(X,T)$ has topological zero entropy, then so does $(\mathcal{M}(X),\widetilde{T})$. We mention that Qiao and Zhou \cite{QiaoZhou17} obtained such a result for the
    notion of sequence entropy.
    
    In another work, Glasner and Weiss \cite{EGlaBWei03} proved that any ergodic system of positive entropy admits {\em every} ergodic system of positive entropy as a quasifactor, which shows, in particular, that the set of quasifactors of an ergodic system of positive entropy is very large.
    
    We also mention that in \cite{Vermersch} the author initiated the investigation on the relationship between the entropy of the MDS $\mathfrak{X}=(X,\mathcal{X},\mu,T)$ and $\widetilde{\mathfrak{X}}=(\mathcal{M}(X),\widetilde{\mathcal{X}},\widetilde{\mu},\widetilde{T})$ in the context of local entropy theory \cite{KerrLi09}. Very recently, Li and Liu, among other findings, expanded it and extended it to amenable group actions \cite{HanfengLi_KairanLiu}.
    \smallskip
    
     %Lindenstrauss and Weiss \cite{Lindenstrauss99, LindenstraussWeiss00} %introduced an important property for a TDS, namely, the Small Boundary %Property (SBP). Let us recall its definition: we say that the TDS $(X,T)$ %has the Small Boundary Property (SBP) if $X$ admits a basis for its %topology consisting of sets $B\in\mathcal{X}$ such that $\nu(\partial %B)=0$ for every $T$-invariant Borel probability measure $\nu$ on $X$. %Clearly, every zero-dimensional TDS has SBP. The significance of the %Small Boundary Property for the present work is that, as shown by %Downarowicz and Huczek in their unifying work about zero-dimensional %extensions \cite{Downarowicz_Huczek}, it follows from deep results due %to Lindenstrauss \cite{Lindenstrauss99} and Lindenstrauss-Weiss %\cite{LindenstraussWeiss00} that, if $(X,T)$ has SBP, then there exists %an \textit{isomorphic} extension $\pi:(Y,S)\rightarrow (X,T)$ from a %zero-dimensional TDS $(Y,S)$ onto $(X,T)$ such that, for every %$T$-invariant Borel probability measure $\mu$ on $X$ there exists a %unique $S$-invariant Borel probability measure $\nu$ on $Y$ such that %$\nu=\mu\circ \pi$. This fact, along with the result on the lifting of %quasifactors due to Glasner-Weiss (Lemma 3.1 from\cite{EGlaBWei03}; %see also Lemma 5.2 from \cite{Meyerovitch11}), allows us to extend our %zero-dimensional results to the SBP setting.
    
    Let $\mathfrak{X}=(X,\mathcal{X},\mu,T)$ be a zero-dimensional MDS and let $A\in\mathcal{X}$ be any clopen set such that $0<\mu(A)<1$. Now, let $0<\eta<1$ and put $\widetilde{A}=\{\nu\in\mathcal{M}(X):\nu(A)>\eta\}$. Since $A\in\mathcal{X}$, $\Psi_{A}(\nu)=\nu(A)$ is a Borel map and $\widetilde{A}=\Psi^{-1}(\eta,1]$, we see that $\widetilde{A}\in\widetilde{\mathcal{X}}$. Write $\mathcal{P}=\{A,A^{c}\}$ and $\widetilde{\mathcal{P}}=\{\widetilde{A},\widetilde{A}^{c}\}$. So, $\mathcal{P}$ is a non-trivial two-set partition of $X$ into Borel sets and, as we shall see, if $\widetilde{\mu}\in Q(\mu)$ has full-support (i.e. if it is positive on the non-empty open sets of $\mathcal{M}(supp\ \mu)$), then $\widetilde{\mathcal{P}}$ is a non-trivial two-set partition of $\mathcal{M}(X)$ into Borel sets (Proposition~\ref{twosetpartition}). Observe that, if $supp\ \mu\neq X$, then it is not possible to find a quasifactor of $\mu$ having full-support on all $\mathcal{M}(X)$. So, when saying that a quasifactor has full-support, it is understood that it has full-support on $\mathcal{M}(supp\ \mu)$. It turns out that our main result (Theorem~\ref{Mainresult}) is based on an analysis of the relationship between the entropy of the stationary stochastic processes generated by the pairs $(\mathfrak{X},\mathcal{P})$ and $(\widetilde{\mathfrak{X}},\widetilde{\mathcal{P}})$, where $\mu$ is ergodic and $\widetilde{\mu}\in Q_{e}(\mu)$ is an ergodic quasifactor of $\mu$ having full-support. In fact, we shall show that, if $\mu$ is ergodic, $\widetilde{\mu}\in Q_{e}(\mu)$ has full-support and $(\widetilde{\mathfrak{X}},\widetilde{\mathcal{P}})$ has zero entropy, then $(\mathfrak{X},\mathcal{P})$ has zero entropy (Corollary~\ref{zeroentropy}). Actually, we show that this fact occurs even {\em continuously} (Theorem~\ref{ergodic_case}). It is worth to point out that our analysis takes advantage of the characterization of the entropy of an ergodic finite-valued stochastic process in terms of the {\em covering-exponent} property. By using that we will show our main result: if  $\mu$ is ergodic with $h_{\mu}(T)>0$, then $h_{\widetilde{\mu}}(\widetilde{T})>0$ for every ergodic quasifactor of $\mu$ having full-support $\widetilde{\mu}\in Q_{e}(\mu)$ (Theorem~\ref{Mainresult}). We remark that we cannot omit the full-support hypothesis for $\widetilde{\mu}\in Q(\mu)$. Actually, we can have $h_{\mu}(T)>0$ and, if we consider $\widetilde{\mu}:=\delta_{\mu}$, then $\widetilde{\mu}\in Q(\mu)$, $\widetilde{\mu}$ is ergodic and $h_{\widetilde{\mu}}(\widetilde{T})=0$. Finally, we conclude this work with a result about the case where $\mu$ and $\widetilde{\mu}$ are not necessarily ergodic. In fact, we show that the present of the process $(\mathfrak{X},\mathcal{P})$ can be predicted with arbitrary accuracy from its past of arbitrarily large finite length, provided that the present of the process $(\widetilde{\mathfrak{X}},\widetilde{\mathcal{P}})$ is sufficiently predictable from its own past of arbitrarily large finite length (Theorem~\ref{approximatedbypast}).

    \section[Preliminaries]{Preliminaries}
    
    Let us recall some definitions and notation from entropy theory.
    In what follows, all logarithms are in base $e$. 
    
    Let $\mathfrak{X}=(X,\mathcal{X},\mu,T)$ be a MDS.
    Given a finite partition $\mathcal{P}=\{P_{0},P_{1},\ldots,P_{k-1}\}$ of $X$, we consider the so-called {\em name map} $\Phi_{\mathcal{P}}:X\to\{0,1,\ldots,k-1\}^{\mathbb{Z}}$ defined by:
    $$(\Phi_{\mathcal{P}}(x))_{n}=j\ \text{if and only if}\  T^{n}x\in P_{j}\ \ (0\leqslant j\leqslant k-1, n\in\mathbb{Z}).
    $$
    The sequence $(\Phi_{\mathcal{P}}(\cdot))_{n\in\mathbb{Z}}$ is a {\em stationary stochastic process}. We say that $(\Phi_{\mathcal{P}}(\cdot))_{n\in\mathbb{Z}}$ is {\em the process generated by $\mathfrak{X}$ and the partition $\mathcal{P}$}.
    
    If $f$ is a random variable in $X$ taking values in $\{0,\ldots,k-1\}$ and we consider, for each $0\leqslant j \leqslant k-1$, the set $P_{j}:=f^{-1}(\{j\})$, then we see that $\mathcal{P}:=\{P_{0},\ldots,P_{k-1}\}$ is a partition of $X$ into Borel sets. Hence, since we can think of a finite partition as a finite-valued random variable that assigns to each point the set containing it, we obtain a correspondence between finite partitions and finite-valued random variables.
    
    The {\em entropy} of a random variable $f$ associated with the finite partition $\mathcal{P}$ is defined by
    
    $$
    H(f):= -\sum_{P\in\mathcal{P}}\mu(P)\log \mu(P),
    $$
    We also write $H(\mathcal{P})=H(f)$.

    Given a stochastic process $(f_{i})_{i\in\mathbb{Z}}$ on $X$ taking values in the finite set $\{0,\ldots,k-1\}$, for each $n\geqslant 1$ we define the {\em joint} of $f_{0},\ldots,f_{n-1}$ by: $$\displaystyle\bigvee_{i=0}^{n-1}f_{i}:=\{P_{0} \cap \cdots \cap P_{n-1}:
    P_{0} \in \mathcal{P}_{0},\ldots,P_{n-1} \in \mathcal{P}_{n-1}\}=\bigvee_{i=0}^{n-1}\mathcal{P}_{i},$$
    where $\mathcal{P}_{i}$ is the partition of $X$ corresponding to $f_{i}$ ($0\leqslant i\leqslant n-1$).

    The {\em entropy of the stochastic process $(f_{i})_{i\in\mathbb{Z}}$} is defined by the following expression:
    $$H((f_{i})_{-\infty}^{+\infty}):= \displaystyle\lim_{n \to \infty}(1/n) H\big(\bigvee_{i=0}^{n-1}\mathcal{P}_{i}\big).
    $$
    
    The {\em entropy of $T$ with respect to $\mathcal{P}$} is defined by
    $$
    h_{\mu}(T,\mathcal{P}):= \displaystyle\lim_{n \to \infty}(1/n) H\big(\bigvee_{i=0}^{n-1}T^{-i}\mathcal{P}\big).
    $$
    Clearly, we have $h_{\mu}(T,\mathcal{P})=H((\Phi_{\mathcal{P}})_{-\infty}^{+\infty})$.
    
    Given $n\in\mathbb{N}$, a finite partition $\mathcal{P}$ of $X$ into Borel sets and $0<\gamma<1$ we denote by $N(n,\mathcal{P},T,\gamma)$ the minimum cardinality of a subcollection $\displaystyle\mathcal{G}\subseteq\bigvee_{i=0}^{n-1}T^{-i}\mathcal{P}$ needed to cover a Borel set $D\subseteq X$ with $\mu(D)\geqslant\gamma$. If $(X,\mathcal{X},\mu,T)$ is \emph{ergodic}, then $h_{\mu}(T,\mathcal{P})$ has the {\em covering-exponent} property; that is,
    
    $$h_{\mu}(T,\mathcal{P})=\displaystyle\lim_{n\to\infty}(1/n)\log N(n,\mathcal{P},T,\gamma)\ \text{for every}\ 0<\gamma<1 $$\ (see, for example, Theorem 5.1 on page 72 from \cite{Rudolph} or Theorem I.7.4 on page 68 from \cite{Shields}).

    Finally, the {\em entropy} of $T$ is given by
    $$
    h_{\mu}(T):= \sup_{\mathcal{P}} h_{\mu}(T,\mathcal{P}),
    $$
    where the supremum is taken over all finite partitions $\mathcal{P}$ of $X$ into Borel sets.
    
    In other words, the entropy of $T$ is the supremum over all the entropies of processes of form $(\Phi_{\mathcal{P}})$ with $\mathcal{P}$ being a finite partition of $X$ into Borel sets.
    
    Let $\Pi\subseteq\mathcal{X}$ be the smallest $\sigma$-algebra containing the collection of all sets $A\in\mathcal{X}$ with $h_{\mu}(T,\{A,A^{c}\})=0$. Pinsker \cite{Pinsker} defined $\Pi$ and showed that:
    
    \begin{itemize}
    	\item[(i)]$T^{-1}\Pi =\Pi$;
    	\item[(ii)]If $\mathcal{F}$ is a $T$-invariant $\sigma$-algebra such that $\mathcal{F}\subseteq\Pi$, then $h_{\mu}(T,\{A,A^{c}\})=0$ for every $A\in\mathcal{F}$.
    \end{itemize}
    Thus, $\Pi$ is the largest $T$-invariant $\sigma$-algebra ``with'' zero entropy. We call $\Pi$ the {\em Pinsker $\sigma$-algebra} of the dynamical system $\mathfrak{X}$. Furthermore, we call the restriction of the dynamical system $\mathfrak{X}$ to $\Pi$ the {\em Pinsker factor} of $\mathfrak{X}$. The Pinsker factor is the {\em deterministic part} of $\mathfrak{X}$. The books by Glasner \cite{EGlabook03} and by Parry \cite{Parry} are standard references for the study of the Pinsker factor.
    Note that $h_{\mu}(T)=0$ if, and only if, $\Pi=\mathcal{X}$. Equivalently, we have $h_{\mu}(T)>0$ if, and only if, there exists some $A\in\mathcal{X}$ with $h_{\mu}(T,\{A,A^{c}\})>0$.
       
    \section{Our results}
    
    We begin with a simple result concerning two-set partitions:
    
    \begin{proposition}{\label{twosetpartition}} Let $(X,\mathcal{X},\mu,T)$ be a zero-dimensional MDS, let $A\in\mathcal{X}$ be any clopen set with $0<\mu(A)<1$, $0<\eta<1$ and let $\widetilde{\mu}\in Q(\mu)$ having full-support. If $\widetilde{A}:=\{\nu\in\mathcal{M}(X):\nu(A)>\eta\}$, then $0<\widetilde{\mu}(\widetilde{A})<1$. 
    	\begin{proof} Suppose $\widetilde{\mu}(\widetilde{A})=0$. So, $\widetilde{\mu}(\{\nu:\nu(A)\leqslant \eta\})=1$. On the other hand, since $\partial A=\varnothing$, it follows that $\widetilde{\mu}(\{\nu:\nu(\overline{A})\leqslant \eta\})=\widetilde{\mu}(\{\nu:\nu(A)\leqslant \eta\})=1$. Since $\widetilde{\mu}(\mathcal{M}(supp\ \mu))=1$, we obtain $\widetilde{\mu}(\{\nu:\nu(\overline{A})\leqslant\eta\}\cap\mathcal{M}(supp\ \mu))=1$. Hence, since $\widetilde{\mu}$ has full-support, it follows that $\{\nu:\nu(\overline{A})\leqslant \eta\}\cap\mathcal{M}(supp\ \mu)$ is dense in $\mathcal{M}(supp\ \mu)$. But, since $\{\nu:\nu(\overline{A})\leqslant \eta\}\cap\mathcal{M}(supp\ \mu)$ is closed in $\mathcal{M}(supp\ \mu)$ (because $A$ is clopen), we obtain $\mathcal{M}(supp\ \mu)=\{\nu:\nu(\overline{A})\leqslant \eta\}\cap\mathcal{M}(supp\ \mu)$, which is impossible. The same argument can be used to show that we cannot have $\widetilde{\mu}(\widetilde{A})=1$ either. Hence, we obtain $0<\widetilde{\mu}(\widetilde{A})<1$. This completes the proof.
    		
    	\end{proof}
    \end{proposition}
    
    Let $(X,\mathcal{X},\mu,T)$ be a zero-dimensional MDS. To avoid unnecessary repetitions, for the rest of the paper we shall consider a clopen set $A\in\mathcal{X}$ with $0<\mu(A)<1$, $0<\eta<1$, $\widetilde{A}=\{\nu\in\mathcal{M}(X):\nu(A)>\eta\}$ and write $\mathcal{P}=\{P_{0},P_{1}\}$ and $\widetilde{\mathcal{P}}=\{\widetilde{P}_{0},\widetilde{P}_{1}\}$, where $P_{0}=A$, $P_{1}=A^{c}=X\setminus A$, $\widetilde{P}_{0}=\widetilde{A}$, $\widetilde{P}_{1}=\widetilde{A}^{c}=\mathcal{M}(X)\setminus\widetilde{A}$. So, $\mathcal{P}:=\{P_{0},P_{1}\}$ is a non-trivial partition of $X$ into Borel sets and, by Proposition~\ref{twosetpartition} we see that $\widetilde{\mathcal{P}}:=\{\widetilde{P}_{0},\widetilde{P}_{1}\}$ is a non-trivial partition of $\mathcal{M}(X)$ into Borel sets.
    
    \begin{lemma}\label{largemeasurelemma}Let $(X,\mathcal{X},\mu,T)$ be a zero-dimensional MDS and let $\widetilde{\mu}\in Q(\mu)$ having full-support. 
    	
    	Given any $0<\gamma<1$ and any $m,n\geqslant 1$ there exists $\gamma^{\prime}=\gamma^{\prime}(\gamma,m,n)>0$  such that, for any functions $\sigma_{1},\ldots,\sigma_{m}\in\{0,1\}^{\{0,\ldots,n-1\}}$, if $\displaystyle\widetilde{\mu}\Big(\bigcup_{l=1}^{m}\bigcap_{i=0}^{n-1}\widetilde{T}^{-i}\widetilde{P}_{\sigma_{l}(i)}\Big)\geqslant \gamma^{\prime}$, then $\displaystyle\mu\Big(\bigcup_{l=1}^{m}\bigcap_{i=0}^{n-1}T^{-i}P_{\sigma_{l}(i)}\Big)\geqslant \gamma$.
    	\begin{proof}To obtain a contradiction, suppose that the result does not hold. In this case, there are  $0<\gamma<1$, $m\geqslant 1$ and $n\geqslant 1$ such that, for every $k\geqslant 1$ there are functions $\sigma_{1}^{k},\ldots,\sigma_{m}^{k}\in\{0,1\}^{\{0,\ldots,n-1\}}$ such that $\displaystyle\widetilde{\mu}\Big(\bigcup_{l=1}^{m}\bigcap_{i=0}^{n-1}\widetilde{T}^{-i}\widetilde{P}_{\sigma_{l}^{k}(i)}\Big)>1- 2^{-k}$ but $\displaystyle\mu\Big(\bigcup_{l=1}^{m}\bigcap_{i=0}^{n-1}T^{-i}P_{\sigma_{l}^{k}(i)}\Big)< \gamma$. Put $\tau_{l}^{k}(i):=\{0,1\}\setminus\{\sigma_{l}^{k}(i)\}$. Since $\displaystyle\sum_{k=1}^{\infty}\displaystyle\widetilde{\mu}\Big(\bigcap_{l=1}^{m}\bigcup_{i=0}^{n-1}\widetilde{T}^{-i}\widetilde{P}_{\tau_{l}^{k}(i)}\Big)< \infty$, it follows from the Borel-Cantelli lemma that , for $\widetilde{\mu}$-a.e. $\nu\in\mathcal{M}(X)$ there exists some $k_{0}=k_{o}(\nu)\geqslant 1$ such that $\nu\notin \displaystyle \bigcap_{l=1}^{m}\bigcup_{i=0}^{n-1}\widetilde{T}^{-i}\widetilde{P}_{\tau_{l}^{k}(i)}$ for all $k\geqslant k_{0}$. In other words, if  
    		$$\widetilde{\mathcal{G}}:=\Big\{\nu\in\mathcal{M}(X):\exists k_{0}\geqslant 1; \nu\in\displaystyle\bigcup_{l=1}^{m}\bigcap_{i=0}^{n-1}\widetilde{T}^{-i}\widetilde{P}_{\sigma_{l}^{k}(i)}\ \text{for all}\ k\geqslant k_{0}\Big\},$$
    		then $\widetilde{\mu}(\widetilde{\mathcal{G}})=1$. 
    		 On the other hand, since $\displaystyle\mu\Big(\bigcap_{l=1}^{m}\bigcup_{i=0}^{n-1}T^{-i}P_{\tau_{l}^{k}(i)}\Big)>1- \gamma$ for all $k\geqslant 1$ we may pick some $x\in supp\ \mu$ such that $x\in\displaystyle\bigcap_{l=1}^{m}\bigcup_{i=0}^{n-1}T^{-i}P_{\tau_{l}^{k}(i)}$ for infinitely many $k's$. Consequently, $\delta_{x}\in\mathcal{M}(supp\ \mu)$ is such that $\delta_{x}\in\displaystyle\bigcap_{l=1}^{m}\bigcup_{i=0}^{n-1}\widetilde{T}^{-i}\widetilde{P}_{\tau_{l}^{k}(i)}$ for infinitely many $k's$. Since $\widetilde{\mu}(\widetilde{\mathcal{G}}\cap\mathcal{M}(supp\ \mu))=1$ and $\widetilde{\mu}$ has full-support, it follows that $\widetilde{\mathcal{G}}\cap\mathcal{M}(supp\ \mu)$ is dense in $\mathcal{M}(supp\ \mu)$. So, for every $\varepsilon>0$ there are $k_{0}\geqslant 1$ and $\nu^{\prime}\in\mathcal{M}(supp\ \mu)$ satisfying:
    		\begin{align}&\nu^{\prime}\in\displaystyle\bigcup_{l=1}^{m}\bigcap_{i=0}^{n-1}\widetilde{T}^{-i}\widetilde{P}_{\sigma_{l}^{k}(i)} \ \text{for all}\ k\geqslant k_{0}\label{nuprimeinLemma}\ \text{and}\\    			
    			&\nu^{\prime}(B)\leqslant \delta_{x}(B^{\varepsilon})+\varepsilon \  \text{for all Borel sets} \ B\subseteq X\label{nuprimedeltaxLemma}.
    		\end{align}
    		
    		Now, fix some $k^{\prime}\geqslant k_{0}$ such that $\delta_{x}\in\displaystyle\bigcap_{l=1}^{m}\bigcup_{i=0}^{n-1}\widetilde{T}^{-i}\widetilde{P}_{\tau_{l}^{k^{\prime}}(i)}$. So, for every $1\leqslant l\leqslant m$ there exists some $0\leqslant i^{\prime}\leqslant n-1$ such that $\widetilde{T}^{i^{\prime}}\delta_{x}\in\widetilde{P}_{\tau_{l}^{k^{\prime}}(i^{\prime})}$. On the other hand, by \eqref{nuprimeinLemma} we see that there exists some $1\leqslant l^{\prime}\leqslant m$  such that $\widetilde{T}^{i}\nu^{\prime}\in\widetilde{P}_{\sigma_{l^{\prime}}^{k^{\prime}}(i)}$ for every $0\leqslant i\leqslant n-1$. In particular, for $l=l^{\prime}$ there exists some $0\leqslant i^{\prime}\leqslant n-1$ such that $\delta_{T^{i^{\prime}}x}(A)>\eta$ and $\nu^{\prime}(T^{-i^{\prime}}A^{c})\geqslant 1-\eta$ if $\tau_{l^{\prime}}^{k^{\prime}}(i^{\prime})=0$, or $\delta_{T^{i^{\prime}}x}(A^{c})\geqslant 1-\eta$ and $\nu^{\prime}(T^{-i^{\prime}}A)>\eta$ if $\tau_{l^{\prime}}^{k^{\prime}}(i^{\prime})=1$. Now, without loss of generality we may assume that $\tau_{l^{\prime}}^{k^{\prime}}(i^{\prime})=0$. In this case, we have $T^{i^{\prime}}x\in A$ and $\nu^{\prime}(T^{-i^{\prime}}A^{c})\geqslant 1-\eta$. Moreover, by \eqref{nuprimedeltaxLemma} with $B=T^{-i^{\prime}}A^{c}$ we get $\delta_{x}((T^{-i^{\prime}}A^{c})^{\varepsilon})\geqslant \nu^{\prime}(T^{-i^{\prime}}A^{c})-\varepsilon\geqslant 1-\eta-\varepsilon>0$, whenever $\varepsilon>0$ is small enough depending on $0<\eta<1$. Therefore, we see that $x\in T^{-i^{\prime}}A$ and $x\in (T^{-i^{\prime}}A^{c})^{\varepsilon}$ for every $\varepsilon>0$ small enough depending on $0<\eta<1$. Since $i^{\prime}=i^{\prime}(\varepsilon)$ and $0\leqslant i^{\prime}(\varepsilon)\leqslant n-1$ for all $\varepsilon>0$,  we may find some sequence $\varepsilon_{j}\to 0$ such that $i^{\prime}(\varepsilon_{j}):=i^{\prime}$ for every $j\geqslant 1$. Then, we get $x\in T^{-i^{\prime}}A$ and $x\in (T^{-i^{\prime}}A^{c})^{\varepsilon_{j}}$ for all $j\geqslant 1$. Therefore, since $T^{-i^{\prime}}A^{c}$ is closed, by letting $j\rightarrow\infty$ we conclude that $x\in T^{-i^{\prime}}A\ \cap (T^{-i^{\prime}}A^{c})=\varnothing$, which is a contradiction.
    			 This proves the lemma.
    	\end{proof}
    \end{lemma}
    
    \smallskip
    
    By analogous reasoning, one can prove the following result:
    
    \begin{lemma}Let $(X,\mathcal{X},\mu,T)$ be a zero-dimensional MDS and let $\widetilde{\mu}\in Q(\mu)$ having full-support. 
    	
    	Given any $0<\gamma<1$ and any $m,n\geqslant 1$ there exists $\gamma^{\prime}=\gamma^{\prime}(\gamma,m,n)>0$  such that, for any functions $\sigma_{1},\ldots,\sigma_{m}\in\{0,1\}^{\{0,\ldots,n-1\}}$, if $\displaystyle\widetilde{\mu}\Big(\bigcup_{l=1}^{m}\bigcap_{i=0}^{n-1}\widetilde{T}^{-i}\widetilde{P}_{\sigma_{l}(i)}\Big)\leqslant \gamma^{\prime}$, then $\displaystyle\mu\Big(\bigcup_{l=1}^{m}\bigcap_{i=0}^{n-1}T^{-i}P_{\sigma_{l}(i)}\Big)\leqslant \gamma$.
    \end{lemma}
    \smallskip
    
    \begin{theorem}\label{coveringtheorem}Let $(X,\mathcal{X},\mu,T)$ be a zero-dimensional MDS and let $\widetilde{\mu}\in Q(\mu)$ having full-support. Given $\alpha>0$ there exists $\beta>0$ with the following property: 
    	
    	Given any $0<\gamma<1$ and any $m,n\geqslant 1$ there exist $\gamma^{\prime}=\gamma^{\prime}(\gamma,m,n)>0$ and $l^{\prime}\geqslant 1$  such that, for every $l\geqslant l^{\prime}$, if $N(l,\widetilde{\mathcal{P}},\widetilde{T},\gamma^{\prime})<e^{l.\beta}$, then $N(l,\mathcal{P},T,\gamma)<e^{l.\alpha}$.
    	\begin{proof}To obtain a contradiction, let us assume that the conclusion does not hold. So, there exists some $\alpha>0$ such that, for every $k\geqslant 1$ there are $0<\gamma_{k}<1$ and two sequences $m_{k},n_{k}\geqslant 1$ such that, for every $\gamma^{\prime}>0$ and every $l^{\prime}\geqslant 1$ we can find $l\geqslant l^{\prime}$ such that $N(l,\widetilde{\mathcal{P}},\widetilde{T},\gamma^{\prime})<e^{l.k^{-1}}$ but $N(l,\mathcal{P},T,\gamma_{k})\geqslant e^{l.\alpha}$. Let $\gamma^{\prime}_{k}>0$ be associated to $\gamma_{k},m_{k}$ and $n_{k}$ according to Lemma~\ref{largemeasurelemma}. Fix any $l^{\prime}\geqslant 1$. So, for some $l\geqslant l^{\prime}$ we have:

    		\begin{equation}\label{ineq.1} N(l,\widetilde{\mathcal{P}},\widetilde{T},\gamma^{\prime}_{k})<e^{l.k^{-1}}
    		\end{equation}
    		and
    		\begin{equation}\label{ineq.2}N(l,\mathcal{P},T,\gamma_{k})\geqslant e^{l.\alpha}.
    		\end{equation}
    		Now, fix $k\geqslant 1$ large enough so that $1/k\leqslant \alpha/2$. Note that (\ref{ineq.1}) means that there exists $\widetilde{D}\subseteq\mathcal{M}(X)$ with $\widetilde{\mu}(\widetilde{D})\geqslant \gamma^{\prime}_{k}$ that admits a collection $\displaystyle\widetilde{\mathcal{G}}\subseteq\bigvee_{i=0}^{l-1}\widetilde{T}^{-i}\widetilde{\mathcal{P}}$ with $|\widetilde{\mathcal{G}}|<e^{l.k^{-1}}$ as a cover $(*)$. Furthermore, note that \eqref{ineq.2} means that for every $D\subseteq X$ with $\mu(D)\geqslant \gamma_{k}$, if $\displaystyle\mathcal{G}\subseteq\bigvee_{i=0}^{l-1}T^{-i}\mathcal{P}$ is a collection that covers $D$, then $|\mathcal{G}|\geqslant e^{l.\alpha}$ $(**)$. Let $\Sigma$ be the collection of functions $\sigma:\{0,\ldots,l-1\}\rightarrow \{0,1\}$ such that, given any $\widetilde{G}\in\widetilde{\mathcal{G}}$ there exists a necessarily unique $\sigma\in\Sigma$ with $\widetilde{G}=\displaystyle\bigcap_{i=0}^{l-1}\widetilde{T}^{-i}\widetilde{P}_{\sigma(i)}$. We now define $\mathcal{G}\subseteq\displaystyle\bigvee_{i=0}^{l-1}T^{-i}\mathcal{P}$ as follows: $G\in\mathcal{G}$ if, and only if, there exists a (necessarily unique) $\sigma\in\Sigma$ such that $G=\displaystyle\bigcap_{i=0}^{l-1}T^{-i}P_{\sigma(i)}$. Put $D:=\displaystyle\bigcup_{G\in\mathcal{G}}G$. Clearly, $\mathcal{G}$ covers $D$ and $|\mathcal{G}|\leqslant|\widetilde{\mathcal{G}}|$. Moreover, since $\displaystyle\widetilde{\mu}\big(\bigcup_{\widetilde{G}\in\widetilde{\mathcal{G}}}\widetilde{G}\big)\geqslant \gamma^{\prime}_{k}$, it follows from Lemma~\ref{largemeasurelemma} that $\mu(D)\geqslant \gamma_{k}$. Therefore, from $(*)$ and $(**)$ we get:
    		$$e^{l.\alpha}\leqslant |\mathcal{G}|\leqslant |\widetilde{\mathcal{G}}|< e^{l.k^{-1}},
    		$$
    		which contradicts the choice $1/k\leqslant \alpha/2$. This concludes the proof of the theorem.
    	\end{proof}
    \end{theorem}
    
    \begin{theorem}\label{ergodic_case}Let $(X,\mathcal{X},\mu,T)$ be a zero-dimensional MDS and let $\widetilde{\mu}\in Q(\mu)$ having full-support. If $\mu$ and $\widetilde{\mu}$ are ergodic, then the following continuity property holds:
    	
    	Given $\alpha>0$ there exists $\beta>0$ such that, if $h_{\widetilde{\mu}}(\widetilde{T},\widetilde{\mathcal{P}})<\beta$, then $h_{\mu}(T,\mathcal{P})<\alpha$.
    	\begin{proof} Let $\alpha>0$ be arbitrary and take $\beta>0$ associated to $\alpha/2$ as in Theorem~\ref{coveringtheorem}. Let $0<\gamma<1$ and $m,n\geqslant 1$ be given and take $\gamma^{\prime}>0$ and $l^{\prime}\geqslant 1$ also as in Theorem~\ref{coveringtheorem}. Suppose $h_{\widetilde{\mu}}(\widetilde{T},\widetilde{\mathcal{P}})<\beta$. There exists $l\geqslant l^{\prime}$ such that $N(l,\widetilde{\mathcal{P}},\widetilde{T},\gamma^{\prime})<e^{l.\beta}$. Hence, by Theorem~\ref{coveringtheorem} we conclude that $N(l,\mathcal{P},T,\gamma)<e^{l.\alpha/2}$. Therefore, we obtain $h_{\mu}(T,\mathcal{P})=\displaystyle\lim_{l\to\infty}(1/l)N(l,\mathcal{P},T,\gamma)\leqslant\alpha/2<\alpha$, as desired.
    	\end{proof}
    	
    \end{theorem}
    
    \medskip
    
    \begin{corollary}\label{zeroentropy}Let $(X,\mathcal{X},\mu,T)$ be an ergodic zero-dimensional MDS and let $\widetilde{\mu}\in Q_{e}(\mu)$ be an ergodic quasifactor of $\mu$ having full-support. If $h_{\widetilde{\mu}}(\widetilde{T},\widetilde{\mathcal{P}})=0$, then $h_{\mu}(T,\mathcal{P})=0$.
    \end{corollary}
    
     \begin{theorem}{\label{Mainresult}} Let $(X,\mathcal{X},\mu,T)$ be a zero-dimensional MDS. If $\mu$ is ergodic and $h_{\mu}(T)>0$, then $h_{\widetilde{\mu}}(\widetilde{T})>0$ for every ergodic quasifactor of $\mu$ having full-support $\widetilde{\mu}\in Q_{e}(\mu)$.
    	
    	\begin{proof} Suppose $h_{\mu}(T)>0$ and let $\widetilde{\mu}\in Q_{e}(\mu)$ be an ergodic quasifactor of $\mu$ having full-support. So, there exists a clopen set $A\subset X$ with $h_{\mu}(T,\mathcal{P})>0$, where $\mathcal{P}=\{A,A^{c}\}$. Fix any $0<\eta<1$ and put $\displaystyle\widetilde{\mathcal{P}}=\{\widetilde{A},\widetilde{A}^{c}\}$, where $\widetilde{A}=\{\nu:\nu(A)>\eta\}$. By Corollary~\ref{zeroentropy} it follows that $h_{\widetilde{\mu}}(\widetilde{T},\widetilde{\mathcal{P}})>0$, which implies $h_{\widetilde{\mu}}(\widetilde{T})>0$, as desired.
    	\end{proof}
    \end{theorem}
    \medskip
    
    Now we turn to the case where both $\mu$ and $\widetilde{\mu}$ are not necessarily ergodic. As we have said before, we shall show that the present of the process $(\mathfrak{X},\mathcal{P})$ can be predicted with arbitrary accuracy from its past of arbitrarily large finite length, provided that the present of the process $(\widetilde{\mathfrak{X}},\widetilde{\mathcal{P}})$ is sufficiently predictable from its own past of arbitrarily large finite length. For this end we introduce the following notation: given two countable partitions $\mathcal{P}$ and $\mathcal{Q}$, given any $m\in\mathbb{N}$ and any $\varepsilon>0$ we write $\mathcal{P}\subseteq_{(m,\varepsilon)}^{\mu}\mathcal{Q}$ to mean that for every $P\in\mathcal{P}$ there exists some union $\bigcup Q$ of $m$ atoms of $\mathcal{Q}$ with $P\subseteq \bigcup Q$ and $\mu(\bigcup Q\setminus P)<\varepsilon$. Finally, we write $\mathcal{P}\subseteq_{(m,0)}^{\mu}\mathcal{Q}$ if, for every $P\in\mathcal{P}$ there exists some union $\bigcup Q$ of $m$ atoms of $\mathcal{Q}$ with $P\subseteq \bigcup Q$ and $\mu(\bigcup Q\setminus P)=0$. Of course, we have $\mathcal{P}\subseteq_{(m,0)}^{\mu}\mathcal{Q}$ if and only if $\mathcal{P}\subseteq_{(m,\varepsilon)}^{\mu}\mathcal{Q}$ for every $\varepsilon>0$ ($m\in\mathbb{N}$).
    
    \medskip
    The ideas behind the proof of our final result parallel with those of Lemma~\ref{largemeasurelemma}; in particular, the use of the classical Borel-Cantelli lemma is crucial again, though the required analysis is more involved.
    
    \begin{theorem}{\label{approximatedbypast}} Let $(X,\mathcal{X},\mu,T)$ be a zero-dimensional MDS and let $\widetilde{\mu}\in Q(\mu)$ having full-support. The following property holds: 
    	
    	Given any $0<\alpha<1$ and any $m,n\geqslant 1$ there exist $\beta=\beta(\alpha,m,n)>0$, $M\geqslant m$ and $N\geqslant n$ such that, for every $m\leqslant m^{\prime}\leqslant M$ and every $n\leqslant n^{\prime}\leqslant N$, if  $\displaystyle\widetilde{\mathcal{P}}\subseteq_{(m^{\prime},\beta)}^{\widetilde{\mu}}\bigvee_{i=1}^{n^{\prime}}\widetilde{T}^{-i}\widetilde{\mathcal{P}}$, then $\displaystyle\mathcal{P}\subseteq_{(m^{\prime},\alpha)}^{\mu}\bigvee_{i=1}^{n^{\prime}}T^{-i}\mathcal{P}$.
    	
    	\begin{proof} Assume that the result does not hold. In this case, there exist $0<\alpha<1$, $m\geqslant 1$ and $n\geqslant 1$ such that, for every $k\geqslant 1$, every $M\geqslant m$ and every $N\geqslant n$ there are $m\leqslant m_{k}^{\prime}\leqslant M$ and $n\leqslant n_{k}^{\prime}\leqslant N$ with
    		\begin{equation}\label{widetildeapproximatedincluded}\displaystyle\widetilde{\mathcal{P}}\subseteq_{(m_{k}^{\prime},2^{-k})}^{\widetilde{\mu}}\bigvee_{i=1}^{n_{k}^{\prime}}\widetilde{T}^{-i}\widetilde{\mathcal{P}}
    		\end{equation}
    		
    		but
    		\begin{equation}\label{notapproximatedincluded}\displaystyle\mathcal{P}\subsetneq_{(m_{k}^{\prime}\alpha)}^{\mu}\bigvee_{i=1}^{n_{k}^{\prime}}T^{-i}\mathcal{P}.
    		\end{equation}
    		Fix any $M\geqslant m$ and any $N\geqslant n$.
    		Observe that (\ref{widetildeapproximatedincluded}) means that:
    		
    		For each $j\in\{0,1\}$ there are functions $\sigma_{j,k}^{1},\ldots,\sigma_{j,k}^{m_{k}^{\prime}}:\{1,\ldots,n_{k}^{\prime}\}\rightarrow\{0,1\}$ such that $\displaystyle\widetilde{P_{j}}\subseteq\bigcup_{l=1}^{m_{k}^{\prime}}\bigcap_{i=1}^{n_{k}^{\prime}}\widetilde{T}^{-i}\widetilde{P}_{\sigma_{j,k}^{l}(i)}$ and $\displaystyle\widetilde{\mu}\Big(\bigcup_{l=1}^{m_{k}^{\prime}}\bigcap_{i=1}^{n_{k}^{\prime}}\widetilde{T}^{-i}\widetilde{P}_{\sigma_{j,k}^{l}(i)}\setminus \widetilde{P}_{j}\Big)<2^{-k}$.
    		
    		Furthermore, observe that (\ref{notapproximatedincluded}) means that: 
    		
    		There exists $j^{\prime}\in\{0,1\}$ such that $\displaystyle\mu\Big(\bigcup_{l=1}^{m_{k}^{\prime}}\bigcap_{i=1}^{n_{k}^{\prime}}T^{-i}P_{\sigma^{l}(i)}\setminus P_{j^{\prime}}\Big)\geqslant\alpha$, whenever $\sigma^{1},\ldots,\sigma^{m_{k}^{\prime}}:\{1,\ldots,n_{k}^{\prime}\}\rightarrow\{0,1\}$ are functions satisfying $\displaystyle P_{j^{\prime}}\subseteq\bigcup_{l=1}^{m_{k}^{\prime}}\bigcap_{i=1}^{n_{k}^{\prime}}T^{-i}P_{\sigma^{l}(i)}$.
    		
    		Fix $j=j^{\prime}$. Since $\displaystyle\sum_{k=1}^{\infty}\widetilde{\mu}\Big(\bigcup_{l=1}^{m_{k}^{\prime}}\bigcap_{i=1}^{n_{k}^{\prime}}\widetilde{T}^{-i}\widetilde{P}_{\sigma_{j^{\prime},k}^{l}(i)}\setminus \widetilde{P}_{j^{\prime}}\Big)<+\infty$, it follows from the Borel-Cantelli lemma that, for $\widetilde{\mu}$-a.e. $\nu\in\mathcal{M}(X)$ there exists some $k_{0}=k_{0}(\nu)\geqslant 1$ such that $\displaystyle\nu\notin\bigcup_{l=1}^{m_{k}^{\prime}}\bigcap_{i=1}^{n_{k}^{\prime}}\widetilde{T}^{-i}\widetilde{P}_{\sigma_{j^{\prime},k}^{l}(i)}\setminus\widetilde{P}_{j^{\prime}}$ for all $k\geqslant k_{0}$. That is, if $\tau_{j^{\prime},k}^{l}(i):=\{0,1\}\setminus \{\sigma_{j^{\prime},k}^{l}(i)\}$ and 
    		$$\widetilde{\mathcal{G}}:=\Big\{\nu\in\mathcal{M}(X):\exists k_{0}\geqslant 1; \nu\in\displaystyle\bigcap_{l=1}^{m_{k}^{\prime}}\bigcup_{i=1}^{n_{k}^{\prime}}\widetilde{T}^{-i}\widetilde{P}_{\tau_{j^{\prime},k}^{l}(i)}\cup P_{j^{\prime}}\ \text{for all}\ k\geqslant k_{0}\Big\},$$
    		then $\widetilde{\mu}(\widetilde{\mathcal{G}})=1$.
    		Let us consider the following sets: $\widetilde{Q}_{0}:=\{\nu:\nu(A)\geqslant\eta\}$ and 
    		$\widetilde{Q}_{1}:=\{\nu:\nu(A^{c})\geqslant1-\eta\}$. Thus, $\widetilde{Q}_{0}$ and $\widetilde{Q}_{1}$ are closed sets with $\widetilde{P}_{0}\subseteq\widetilde{Q}_{0}$ and $\widetilde{P}_{1}\subseteq\widetilde{Q}_{1}$. Let us consider the following set:
    		$$\widetilde{\mathcal{H}}:=\Big\{\nu\in\mathcal{M}(X):\exists k_{0}\geqslant 1; \nu\in\displaystyle\bigcap_{l=1}^{m_{k}^{\prime}}\bigcup_{i=1}^{n_{k}^{\prime}}\widetilde{T}^{-i}\widetilde{Q}_{\tau_{j^{\prime},k}^{l}(i)}\cup Q_{j^{\prime}}\ \text{for all}\ k\geqslant k_{0}\Big\}.$$
    		Clearly, $\widetilde{\mathcal{G}}\subseteq\widetilde{\mathcal{H}}$ and so, $\widetilde{\mu}(\widetilde{\mathcal{H}})=1$. Hence, $\widetilde{\mu}(\widetilde{\mathcal{H}}\cap\mathcal{M}(supp\ \mu))=1$.  Since $\widetilde{\mu}$ has full-support, it follows that $\widetilde{\mathcal{H}}\cap\mathcal{M}(supp\ \mu)$ is dense in $\mathcal{M}(supp\ \mu)$. Now, let us consider the set 
    		$$\widetilde{\Lambda}:=\Big\{\nu\in\mathcal{M}(X):\exists k_{0}\geqslant 1; \nu\in\displaystyle\bigcap_{l=1}^{m_{k}^{\prime}}\bigcup_{i=1}^{n_{k}^{\prime}}\widetilde{T}^{-i}\widetilde{Q}_{\tau_{j^{\prime},k}^{l}(i)}\ \text{for all}\ k\geqslant k_{0}\Big\}.$$
    		Clearly, $\widetilde{\mathcal{H}}=\widetilde{\Lambda}\cup\widetilde{Q}_{j^{\prime}}$ and so, $\widetilde{\mathcal{H}}\setminus\widetilde{Q}_{j^{\prime}}\subseteq\widetilde{\Lambda}$. Since $\widetilde{Q}_{j^{\prime}}$ is closed, we see that every $\nu\in\mathcal{M}(supp\ \mu)\setminus\widetilde{Q}_{j^{\prime}}$ can be arbitrarily well-approximated by elements from $\widetilde{\Lambda}\cap\mathcal{M}(supp\ \mu)$. 
    		
    		More precisely, we have that for every $\nu\in\mathcal{M}(supp\ \mu)\setminus\widetilde{Q}_{j^{\prime}}$ and every $\varepsilon>0$ there are $k_{0}\geqslant 1$ and $\nu^{\prime}\in\mathcal{M}(supp\ \mu)$ satisfying:
    		\begin{align*}&\nu^{\prime}\in\displaystyle\bigcap_{l=1}^{m_{k}^{\prime}}\bigcup_{i=1}^{n_{k}^{\prime}}\widetilde{T}^{-i}\widetilde{Q}_{\tau_{j^{\prime},k}^{l}(i)} \ \text{for all}\
    			k\geqslant k_{0}\ \text{and}\\
    			&\nu^{\prime}(B)\leqslant \nu(B^{\varepsilon})+\varepsilon \  \text{for all Borel sets} \ B\subseteq X.
    		\end{align*}
    		Without loss of generality we may assume that $j^{\prime}=0$. Put $j^{\prime\prime}:=\{0,1\}\setminus\{j^{\prime}\}$; so, $j^{\prime\prime}=1$. Thus, we can rewrite the above condition as follows:
    		
    		$(*)$ For every $\nu\in\mathcal{M}(supp\ \mu)$ such that $\nu(A)<\eta$ and every $\varepsilon>0$ there are $k_{0}\geqslant 1$ and  $\nu^{\prime}\in\mathcal{M}(supp\ \mu)$ satisfying:
    		\begin{align*}&\nu^{\prime}\in\displaystyle\bigcap_{l=1}^{m_{k}^{\prime}}\bigcup_{i=1}^{n_{k}^{\prime}}\widetilde{T}^{-i}\widetilde{Q}_{\tau_{j^{\prime},k}^{l}(i)} \ \text{for all}\ k\geqslant k_{0}\ \text{and}\\
    			&\nu^{\prime}(B)\leqslant \nu(B^{\varepsilon})+\varepsilon \  \text{for all Borel sets} \ B\subseteq X.
    		\end{align*}
    		On the other hand, since $\displaystyle\widetilde{P}_ {j^{\prime}}\subseteq\bigcup_{l=1}^{m_{k}^{\prime}}\bigcap_{i=1}^{n_{k}^{\prime}}\widetilde{T}^{-i}\widetilde{P}_{\sigma_{j^{\prime},k}^{l}(i)}$ implies $\displaystyle P_{j^{\prime}}\subseteq\bigcup_{l=1}^{m_{k}^{\prime}}\bigcap_{i=1}^{n_{k}^{\prime}}T^{-i}P_{\sigma_{j^{\prime},k}^{l}(i)}$, it follows from (\ref{notapproximatedincluded}) that $\displaystyle\mu\big(\bigcup_{l=1}^{m_{k}^{\prime}}\bigcap_{i=1}^{n_{k}^{\prime}}T^{-i}P_{\sigma_{j^{\prime},k}^{l}(i)}\setminus P_{j^{\prime}}\big)\geqslant\alpha$ for every $k\geqslant 1$. Therefore, we may pick some $x\in supp\ \mu$ such that $x\in\displaystyle\bigcup_{l=1}^{m_{k}^{\prime}}\bigcap_{i=1}^{n_{k}^{\prime}}T^{-i}P_{\sigma_{j^{\prime},k}^{l}(i)}\setminus P_{j^{\prime}}$ for infinitely many $k's$. Consequently, $\delta_{x}\in\mathcal{M}(supp\ \mu)$ is such that $\delta_{x}\in\displaystyle\bigcup_{l=1}^{m_{k}^{\prime}}\bigcap_{i=1}^{n_{k}^{\prime}}\widetilde{T}^{-i}\widetilde{P}_{\sigma_{j^{\prime},k}^{l}(i)}\cap \widetilde{P}_{j^{\prime\prime}}$ for infinitely many $k's$. Furthermore, since $x\in P_{1}=A^{c}$ we have $\delta_{x}\in\mathcal{M}(supp\ \mu)\setminus\widetilde{Q}_{0}$. Therefore, by $(*)$ above with $\nu=\delta_{x}$ we see that for every $\varepsilon>0$ there are $k_{0}\geqslant 1$ and some $\nu^{\prime}\in\mathcal{M}(supp\ \mu)$ such that:
    		\begin{align}&\nu^{\prime}\in\displaystyle\bigcap_{l=1}^{m_{k}^{\prime}}\bigcup_{i=1}^{n_{k}^{\prime}}\widetilde{T}^{-i}\widetilde{Q}_{\tau_{j^{\prime},k}^{l}(i)} \ \text{for all}\ k\geqslant k_{0}\label{nuprimein}\ \text{and}\\
    			&\nu^{\prime}(B)\leqslant \delta_{x}(B^{\varepsilon})+\varepsilon \  \text{for all Borel sets} \ B\subseteq X\label{nuprimedeltax}.
    		\end{align}
    		Now, fix some $k^{\prime}\geqslant k_{0}$ such that $\delta_{x}\in\displaystyle\bigcup_{l=1}^{m^{\prime}_{k^{\prime}}}\bigcap_{i=1}^{n^{\prime}_{k^{\prime}}}\widetilde{T}^{-i}\widetilde{P}_{\sigma_{j^{\prime},k^{\prime}}^{l}(i)}\cap \widetilde{P}_{1}$. So, there exists some $1\leqslant l^{\prime}\leqslant m^{\prime}_{k^{\prime}}$ such that $\widetilde{T}^{i}\delta_{x}\in\widetilde{P}_{\sigma_{j^{\prime},k^{\prime}}^{l^{\prime}}(i)}$ for every $1\leqslant i\leqslant n^{\prime}_{k^{\prime}}$. That is, we have $\delta_{T^{i}x}(A)>\eta$ if $\sigma_{j^{\prime},k^{\prime}}^{l^{\prime}}(i)=0$ and $\delta_{T^{i}x}(A^{c})\geqslant 1-\eta$ if $\sigma_{j^{\prime},k^{\prime}}^{l^{\prime}}(i)=1$ ($1\leqslant i\leqslant n^{\prime}_{k^{\prime}}$). On the other hand, by \eqref{nuprimein} we see that for every $1\leqslant l\leqslant m^{\prime}_{k^{\prime}}$ there exists some $1\leqslant i^{\prime}\leqslant n^{\prime}_{k^{\prime}}$ such that $\widetilde{T}^{i^{\prime}}\nu^{\prime}\in\widetilde{Q}_{\tau_{j^{\prime},k^{\prime}}^{l}(i^{\prime})}$. In particular, for $l=l^{\prime}$ there exists some $1\leqslant i^{\prime}\leqslant n^{\prime}_{k^{\prime}}$ such that $\nu^{\prime}(T^{-i^{\prime}}A)\geqslant \eta$ if $\tau_{j^{\prime},k^{\prime}}^{l^{\prime}}(i^{\prime})=0$ or $\nu^{\prime}(T^{-i^{\prime}}A^{c})\geqslant 1-\eta$ if $\tau_{j^{\prime},k^{\prime}}^{l^{\prime}}(i^{\prime})=1$. Now, without loss of generality we may assume that $\sigma_{j^{\prime},k^{\prime}}^{l^{\prime}}(i^{\prime})=0$ (which is the same as $\tau_{j^{\prime},k^{\prime}}^{l^{\prime}}(i^{\prime})=1)$. In this case, we have $T^{i^{\prime}}x\in A$ and $\nu^{\prime}(T^{-i^{\prime}}A^{c})\geqslant 1-\eta$. Moreover, by \eqref{nuprimedeltax} with $B=T^{-i^{\prime}}A^{c}$ we get $\delta_{x}((T^{-i^{\prime}}A^{c})^{\varepsilon})\geqslant \nu^{\prime}(T^{-i^{\prime}}A^{c})-\varepsilon\geqslant 1-\eta-\varepsilon>0$, whenever $\varepsilon>0$ is small enough depending on $0<\eta<1$. Therefore, we see that $x\in T^{-i^{\prime}}A$ and $x\in (T^{-i^{\prime}}A^{c})^{\varepsilon}$ for every $\varepsilon>0$ small enough depending on $0<\eta<1$. Since $i^{\prime}=i^{\prime}(\varepsilon)$ and $1\leqslant i^{\prime}(\varepsilon)\leqslant N$ for all $\varepsilon>0$,  we may find some sequence $\varepsilon_{j}\to 0$ such that $i^{\prime}(\varepsilon_{j}):=i^{\prime}$ for every $j\geqslant 1$. Then, we get $x\in T^{-i^{\prime}}A$ and $x\in (T^{-i^{\prime}}A^{c})^{\varepsilon_{j}}$ for all $j\geqslant 1$. Therefore, since $T^{-i^{\prime}}A^{c}$ is closed, by letting $j\rightarrow\infty$ we conclude that $x\in T^{-i^{\prime}}A\ \cap (T^{-i^{\prime}}A^{c})=\varnothing$, which is a contradiction.
    		This concludes the proof of the theorem.
    		\end{proof}
    	\end{theorem}

    \smallskip

    \section*{Acknowledgments}
    	The author would like to dedicate this work to his daughter, Laura C. Vermersch, in gratitude for the joyful moments during its writing. %The author would like to thank the anonymous referee for helpful comments that improved the manuscript.

   % \end{Backmatter}
    
    \end{document}